# An Existence Theorem for a Class of Wrinkling Models for Highly Stretched Elastic Sheets

Timothy J. Healey

Department of Mathematics, Cornell University, Ithaca, NY 14853

tjh10@cornell.edu

**Abstract**

We consider a class of models motivated by previous numerical studies of wrinkling in highly stretched, thin rectangular elastomer sheets. The model is characterized by a finite-strain hyperelastic membrane energy perturbed by small bending energy. In the absence of the latter, the membrane energy density is not rank-one convex for general spatial deformations but reduces to a polyconvex function when restricted to planar deformations, i.e., two-dimensional hyperelasticity. In addition, it grows unbounded as the local area ratio approaches zero. The small-bending component of the model is the same as that in the classical von Kármán model. The latter penalizes arbitrarily fine-scale wrinkling, resolving both the amplitude and wavelength of wrinkles. Here we prove the existence of energy minima for a general class of such models.

## 1. Introduction

We consider a class of models characterized by properties first proposed in [8] for the study of wrinkling in highly stretched, thin rectangular elastomer sheets. The formally derived model comprises a finite-strain hyperelastic membrane perturbed by small bending energy. The latter penalizes arbitrarily fine-scale wrinkling, resolving both the amplitude and wavelength of wrinkles. Numerical bifurcation-continuation methods are employed in the analysis of a finite-element discretization. The approach successfully predicts stable families of equilibria, and new phenomena are uncovered. However, there is no apparent way to realize those results rigorously via bifurcation theory. For instance, even the determination of the planar, unwrinkled solution family presents an unsolved problem of two-dimensional (2D) nonlinear elasticity. Clearly, the analysis of bifurcation to non-planar wrinkled states is ruled out. In addition, a rigorous derivation of the model via $\Gamma$- convergence in the limit of zero thickness seems unlikely. We also mention that the model proposed in [8] provides accurate predictions in comparison with experiment within the context of purely elastic behavior, cf. [5].

Here we prove the existence of energy minimizing configurations for a general class of models motivated by that employed in [8]. Compelling numerical evidence aside, we subscribe here to the philosophy espoused in [3][1]. The membrane energy employed in [8] is the same as that used by Treloar to interpret bi-axial experiments on thin rubber sheets [15], using Mooney's correction to the neo-Hookean model [10], cf. [11, Eq. (3.22)]. In the absence of bending energy, it is not rank-one convex but reduces to a polyconvex stored energy when restricted to the plane, i.e., 2D hyperelasticity. (See [1] for the terminology.) In addition, it grows unbounded as the local area ratio approaches zero. The additive small-bending component of the model used in [8] is the same as that in the classical von Kármán model [16], presuming small out-of-plane displacements and slopes.

---

[1] "Only a mathematical existence proof can ensure that the mathematical description of a physical phenomenon is meaningful".

Several different models and approaches to wrinkling have been proposed in the literature. We mention some that incorporate finite elasticity in the membrane part of the model, e.g., [2], [6], [7], [12], and [14]. With the possible exception of [12], each of these has the same basic structure as that in [8], viz., a hyperelastic membrane plus small bending energy. On the other hand, most of them feature properties different from that used in [8]. A geometrically exact St.Venant-Kirchhoff membrane energy is employed in both [7] and [14]. It remains bounded as the local area ratio approaches zero; even when restricted to planar behavior, the membrane energy is not rank-one convex. The membrane model in [2] is derived from 3D hyperelasticity presuming isotropy and convexity in the right stretch tensor. The shortcomings of that convexity assumption within the context of rubber elasticity are well known [1], and the resulting membrane energy density has the same properties as the St. Venant-Kirchhoff membrane model just described. The model used in [14] is formally derived in [13] from 3D hyperelasticity via small-thickness asymptotics. It incorporates the exact relative-curvature tensor in the quadratic bending energy, whereas the approximate von Kármán bending model is employed in [2], [6], [7] and [8]. The model in [6] is the same as that proposed in [8]; the effects of various material properties on wrinkling behavior are explored in the former. A finite-deformation, hyperelastic thin-shell model from a commercial finite-element code is employed in [12]. Except for [2], all the above-mentioned wrinkling studies are numerical. In all cases, the correct qualitative phenomenon is apparently captured. However, the existence of solutions is not addressed in any of these.

The outline of the work is as follows. In Section 2 we present the formulation for the class of problems to be treated. We first recall the geometric approximation associated with small bending. The same assumption impacts the measure of the local-area ratio, which is consistently approximated by its value based solely on the in-plane component of the deformation. Both approximations are incorporated in the model. We then present our constitutive hypotheses, generalizing the properties of the model from [8] as discussed above. Our main result is presented in Section 3, establishing the existence of an energy-minimizing configuration maintaining local orientation. The latter is based on the assumed blow up of the energy density as the approximate local area ratio just described approaches zero. We make some final remarks in Section 4. In particular, we observe that a class of partial first variations in the out-of-plane direction can be taken rigorously at an energy minimizer. We also observe that our approach to existence does not extend to any of the other models discussed above.

## 2. Formulation

We follow the common practice of letting $\mathbb{R}^3$ denote both Euclidean point space and its translate inner-product space. The Euclidean inner-product is denoted by $\mathbf{a} \cdot \mathbf{b}$ for all $\mathbf{a}, \mathbf{b} \in \mathbb{R}^3$. We let $\{\mathbf{e}_1, \mathbf{e}_2, \mathbf{e}_3\}$ denote the standard, right-handed orthonormal basis for $\mathbb{R}^3$, and we henceforth make the identification $\mathbb{R}^2 \cong span\{\mathbf{e}_1, \mathbf{e}_2\}$. Let $\Omega \subset \mathbb{R}^2$ denote a bounded domain with a strongly Lipschitz boundary $\partial \Omega$; we identify $\bar{\Omega}$ with the reference configuration for a material surface in a flat state. We decompose a deformation $\mathbf{f} : \bar{\Omega} \to \mathbb{R}^3$ according to

$$\mathbf{f}(\mathbf{x}) = \mathbf{h}(\mathbf{x}) + w(\mathbf{x})\mathbf{e}_3, \tag{1}$$

where $\mathbf{h} : \bar{\Omega} \to \mathbb{R}^2$. The deformation gradient of $\mathbf{f}$ is then given by

$$\mathbf{F} := \nabla \mathbf{f} = \mathbf{f}_{,\alpha} \otimes \mathbf{e}_\alpha = \nabla \mathbf{h} + \mathbf{e}_3 \otimes \nabla w, \tag{2}$$

where $\mathbf{H} := \nabla \mathbf{h} = \mathbf{h}_{,\alpha} \otimes \mathbf{e}_\alpha$ (summation convention with Greek indices summing from 1-2) and $\nabla w$ denotes the usual gradient of a scalar-valued function. The second gradient of $\mathbf{f}$ is given by

$$\mathbf{G} := \nabla^2 \mathbf{f} = \mathbf{f},_{\alpha\beta} \otimes \mathbf{e}_\alpha \otimes \mathbf{e}_\beta = \nabla^2 \mathbf{h} + \mathbf{e}_3 \otimes \nabla^2 w, \tag{3}$$

where $\nabla^2 \mathbf{h} = \mathbf{h},_{\alpha\beta} \otimes \mathbf{e}_\alpha \otimes \mathbf{e}_\beta$ and $\nabla^2 w = w,_{\alpha\beta} \mathbf{e}_\alpha \otimes \mathbf{e}_\beta$. The triple tensor product employed in (3) is defined by

$$(\mathbf{a} \otimes \mathbf{b} \otimes \mathbf{c})\mathbf{u} := \mathbf{a} \otimes \mathbf{b}(\mathbf{c} \cdot \mathbf{u}),$$
$$(\mathbf{a} \otimes \mathbf{b} \otimes \mathbf{c})(\mathbf{u} \otimes \mathbf{v}) := \mathbf{a}(\mathbf{b} \cdot \mathbf{u})(\mathbf{c} \cdot \mathbf{v}),$$
$$\mathbf{u} \cdot (\mathbf{a} \otimes \mathbf{b} \otimes \mathbf{c}) := (\mathbf{u} \cdot \mathbf{a})\mathbf{b} \otimes \mathbf{c},$$

for all $\mathbf{a}, \mathbf{b}, \mathbf{c}, \mathbf{u}, \mathbf{v} \in \mathbb{R}^3$.

Let $\mathbf{n}$ denote the unit-normal field on the deformed surface $\mathbf{f}(\Omega)$ that points in the same direction as the vector field $\mathbf{f},_1 \times \mathbf{f},_2$. The local area ratio is then given by

$$J := \sqrt{\det[\mathbf{F}^T \mathbf{F}]} = \mathbf{n} \cdot (\mathbf{f},_1 \times \mathbf{f},_2). \tag{4}$$

The relative curvature tensor is defined by

$$\mathbf{K} := -\mathbf{n} \cdot \nabla^2 \mathbf{f} = -(\mathbf{n} \cdot \mathbf{f},_{\alpha\beta})\mathbf{e}_\alpha \otimes \mathbf{e}_\beta. \tag{5}$$

As in [8], we make a small-curvature (von Kármán) approximation, via $\mathbf{n} \simeq \mathbf{e}_3$, in which case (5) via (3) yields

$$\mathbf{K} \simeq \nabla^2 w, \tag{6}$$

where we have used the fact that $\mathbf{e}_3 \cdot \mathbf{h},_{\alpha\beta} \equiv 0$. For consistency, (4) then becomes

$$J \simeq \mathbf{e}_3 \cdot (\mathbf{f},_1 \times \mathbf{f},_2) = \det \nabla \mathbf{h}. \tag{7}$$

*We henceforth adopt* (6) *and* (7) *as definitions*.

We assume the existence of an energy density function

$$W(\mathbf{K}, \mathbf{F}), \quad W : S^{2\times 2} \times L_+^{3\times 2} \to [0, \infty), \tag{8}$$

where $S^{2\times 2} := \{\mathbf{A} \in \mathbb{R}^{2\times 2} : \mathbf{A}^T = \mathbf{A}\}$ and $L_+^{3\times 2} := \{\mathbf{F} \in \mathbb{R}^{3\times 2} : J > 0\}$. To motivate our forthcoming assumptions, we recall the energy density employed in [8], which combines finite nonlinear elasticity for the membrane part plus small bending energy:

$$W(\mathbf{K}, \mathbf{F}) = \Psi(\mathbf{K}) + \Phi(\mathbf{C}), \tag{9}$$

where $\mathbf{C} = \mathbf{F}^T \mathbf{F}$ is the $2 \times 2$ right Cauchy-Green strain tensor. The bending energy $\Psi$ is the same as that in the von Kármán theory [16], while the membrane energy $\Phi$ is based on the Mooney-Rivlin model for a thin membrane, cf. [11]. In particular, $\mathbb{R}^{3\times 2} \times (0, \infty) \ni (\mathbf{F}, J) \mapsto \Phi$ is convex, with $J$ as in (4). In addition, the property

$$\Phi(\mathbf{C}) \nearrow \infty \text{ as } \det \mathbf{C} \searrow 0, \tag{10}$$

arises naturally. Although it is polyconvex when $\mathbf{F}$ is restricted to $\mathbb{R}^{2\times 2}$, i.e., 2D hyperelasticity, $\Phi$ is not rank-one convex for spatial deformations, i.e., $\mathbf{F} \in \mathbb{R}^{3\times 2}$.

We make the following constitutive hypotheses (based on (6) and (7)):

(H1) There is a convex $C^1$ function $\Gamma : \mathbb{R}^{2\times 2} \times \mathbb{R}^{3\times 2} \times (0,\infty) \to [0,\infty)$ such that

$$W(\mathbf{K},\mathbf{F}) \equiv \Gamma(\mathbf{K},\mathbf{F},J).$$

(H2) There exist constants $p>1$, $q>4/3$, $p>2q/(2+q)$ for $q \geq 2$, $r>1$, $C_1>0$ and $C_2$ such that

$$W(\mathbf{K},\mathbf{F}) \geq C_1 \left\{ |\mathbf{K}|^p + |\mathbf{F}|^q + J^r \right\} + C_2.$$

(H3) $\Gamma \nearrow \infty$ as $J \searrow 0$.

We let $L^p(\Omega,\mathbb{R}^n)$ and $W^{k,p}(\Omega,\mathbb{R}^n)$ denote the usual Lebesgue and Sobolev spaces, respectively, each equipped with the standard norm. The energy functional is defined as

$$E[\mathbf{h},w] := \int_\Omega W(\nabla^2 w, \nabla \mathbf{h} + \mathbf{e}_3 \otimes \nabla w)dx - \varphi[w,\mathbf{h}], \qquad (11)$$

where $\varphi$ is a bounded linear functional on $W^{2,p}(\Omega) \times W^{1,q}(\Omega,\mathbb{R}^2)$, representing dead loading. For example,

$$\varphi[\mathbf{h},w] = \int_\Omega [\mathbf{m}\cdot \nabla w + \mathbf{b}\cdot(\mathbf{h}+w\mathbf{e}_3)]dx,$$

where $\mathbf{m},\mathbf{b}\in L^\infty(\Omega,\mathbb{R}^2)$ are prescribed.

The set of admissible configurations is defined as

$$\begin{aligned}
\mathcal{A} := \{(\mathbf{h},w) &\in W^{1,q}(\Omega,\mathbb{R}^2)\times W^{2,p}(\Omega): J\in L^r(\Omega),\ J>0 \\
&\text{a.e. in } \Omega,\ \mathbf{h}-\mathbf{h}_o \in W^{1,q}_\Gamma(\Omega,\mathbb{R}^2),\ w-w_o \in W^{2,p}_\Gamma(\Omega)\},
\end{aligned} \qquad (12)$$

where $(\mathbf{h}_o,w_o)\in W^{1,q}(\Omega,\mathbb{R}^2)\times W^{2,p}(\Omega)$ is prescribed. Here we define

$$W^{2,p}_\Gamma(\Omega) := \{u\in W^{2,p}(\Omega): u = \nabla u \cdot \mathbf{v} = 0 \text{ on } \Gamma\},$$
$$W^{1,q}_\Gamma(\Omega,\mathbb{R}^2) := \{\mathbf{g}\in W^{1,q}(\Omega,\mathbb{R}^2): \mathbf{g} = \mathbf{0} \text{ on } \Gamma\},$$

where $\Gamma \subset \partial\Omega$ has surface measure $|\Gamma|>0$ and $\mathbf{v}$ denotes the outward unit-normal field on $\Gamma$. The boundary prescriptions above are understood in the sense of trace.

**Remark 1** Although it plays no direct role in what follows, we claim that material objectivity here reads $W(\mathbf{K},\mathbf{QF}) = W(\mathbf{K},\mathbf{F})$ $\mathbf{Q}\in SO(3)$. To see this, note that the exact relative curvature tensor (5) is invariant under rotations: $\mathbf{n}\to \mathbf{Qn}$, $\mathbf{G}\to \mathbf{QG} \Rightarrow \mathbf{Qn}\cdot \mathbf{QG} = \mathbf{n}\cdot \mathbf{G}$ $\forall\ \mathbf{Q}\in SO(3)$, which also holds for the special case $\mathbf{n}=\mathbf{e}_3$.

## 3. Existence of Energy Minimizers

Our main result is the following:

**Theorem 1** *Assume that $(\mathbf{h}_o,w_o)\in \mathcal{A}$ and $E[\mathbf{h}_o,w_o]<\infty$. Then there exists $(\mathbf{h}_*,w_*)\in \mathcal{A}$ such that $E[\mathbf{h}_*,w_*] = \inf_\mathcal{A} E[\mathbf{h},w] < \infty.$*

***Proof*** From (H2) and embedding, it follows that $w\in W^{2,p}(\Omega)\cap W^{1,q}(\Omega)$, and integration of the growth condition implies

$$\int_\Omega W(\nabla^2 w, \nabla \mathbf{h} + \mathbf{e}_3 \otimes \nabla w) dx$$
$$\geq C_3 \left\{ \|\nabla^2 w\|_{L^p}^p + \|\nabla \mathbf{h} + \mathbf{e}_3 \otimes \nabla w\|_{L^q}^q + \|J\|_{L^r}^r \right\} + C_4,$$

for constants $C_3 > 0$, $C_4$. Accounting for the loading term $\varphi[\cdot]$ and employing Poincaré's inequality, we then arrive at

$$E[w, \mathbf{h}] \geq C_5 \left\{ \|w\|_{W^{2,p}}^p + \|(\mathbf{h}, w)\|_{W^{1,q}}^q + \|J\|_{L^r}^r \right\} - C_6 \left\{ \|w\|_{W^{2,p}} + \|\mathbf{h}\|_{W^{1,q}} \right\} + C_7, \tag{13}$$

for some constants $C_5, C_6 > 0$, and $C_7$. Since $p, q > 1$, (13) implies there are constants $M > 0$ and $D$ that

$$E[w, \mathbf{h}] \geq M \left\{ \|w\|_{W^{2,p}}^p + \|(\mathbf{h}, w)\|_{W^{1,q}}^q + \|J\|_{L^r}^r \right\} + D. \tag{14}$$

Next, let $\{(\mathbf{h}_k, w_k)\} \subset \mathcal{A}$ be a minimizing sequence for $E[\cdot]$, i.e.,

$$\lim_{k \to \infty} E[\mathbf{h}_k, w_k] = \inf_{\mathcal{A}} E[\mathbf{h}, w].$$

In consonance with (7), we define $J_k := \det \nabla \mathbf{h}_k$. Given $\inf_{\mathcal{A}} E[w, \mathbf{h}] < \infty$, the lower bound (14) implies that the sequences $\{w_k\}, \{(\mathbf{h}_k, w_k)\}$ and $\{J_k\}$ are uniformly bounded in the reflexive Banach spaces $W^{2,p}(\Omega)$, $W^{1,q}(\Omega, \mathbb{R}^3)$ and $L^r(\Omega)$, respectively. Hence, there are weakly convergent subsequences (not relabeled) such that $w_k \rightharpoonup w_* \in W^{2,p}(\Omega)$, $(\mathbf{h}_k, w_k) \rightharpoonup (\mathbf{h}_*, w_*) \in W^{1,q}(\Omega, \mathbb{R}^3)$, and $J_k \rightharpoonup \alpha_* \in L^r(\Omega)$. For $\mathbf{h} \in W^{1,q}(\Omega, \mathbb{R}^2)$ with $q \geq 2$, the determinant (7) is a well-defined element of $L^{q/2}(\Omega)$. Otherwise, for $4/3 < q < 2$, the distributional determinant is defined via

$$\int_\Omega J_k \phi \, dx := -\frac{1}{2} \int_\Omega ([Cof \nabla \mathbf{h}_k]^T \mathbf{h}_k) \cdot \nabla \phi \, dx, \quad \forall \phi \in C_c^\infty(\Omega), \tag{15}$$

where $Cof\, \mathbf{A}$ denotes the cofactor matrix of the $2 \times 2$ matrix $\mathbf{A}$. With $q > 4/3$, it is well known that the right side of (15) converges in the sense of distributions [4], i.e.,

$$\int_\Omega J_k \phi \, dx \to \int_\Omega J_* \phi \, dx, \quad \forall \phi \in C_c^\infty(\Omega), \tag{16}$$

where $J_* = \det \nabla \mathbf{h}_*$ (interpreted according to (15)). Comparing the weak convergence $J_k \rightharpoonup \alpha_* \in L^r(\Omega)$ to (16), we conclude that $\alpha_* = J_* \in L^r(\Omega)$.

To show that $(\mathbf{h}_*, w_*) \in \mathcal{A}$, we first claim that $J_* > 0$ a.e. in $\Omega$. Indeed, we may construct a sequence of convex combinations of $\{J_k\}$ that converges strongly to $J_*$ in $L^1(\Omega)$ and thus converges (possibly as a subsequence) to $J_*$ a.e. in $\Omega$. Now each $J_k > 0 \Rightarrow J_* \geq 0$ a.e. in $\Omega$. Then (H3) and Fatou's lemma contradicts the possibility of $J_* = 0$ a.e. in $\Lambda \subset \Omega$ with $|\Lambda| > 0$. Finally, the closed subspaces $W_\Gamma^{2,p}(\Omega)$ and $W_\Gamma^{1,q}(\Omega, \mathbb{R}^2)$ are weakly closed, implying that $w_* - w_o \in W_\Gamma^{2,p}(\Omega)$ and $\mathbf{h}_* - \mathbf{h}_o \in W_\Gamma^{1,q}(\Omega, \mathbb{R}^2)$.

To finish, we need to ensure that $E[\cdot]$ is weakly lower semi-continuous, viz.,

$$\liminf_{j \to \infty} E[\mathbf{h}_j, w_j] \geq E[\mathbf{h}, w], \tag{17}$$

whenever $w_j \rightharpoonup w$ in $W^{2,p}(\Omega)$, $(\mathbf{h}_k, w_k) \rightharpoonup (\mathbf{h}, w)$ in $W^{1,q}(\Omega, \mathbb{R}^2)$, and $J_k \rightharpoonup J$ in $L^r(\Omega)$, for $J_k, J > 0$ a.e. in $\Omega$. Since, $\varphi[\cdot]$ is a bounded linear functional, it's enough to focus on the weak lower semi-continuity of the internal energy functional

$$I[\mathbf{h}, w] := \int_\Omega W(\nabla^2 w, \nabla \mathbf{h} + \mathbf{e}_3 \otimes \nabla w) dx.$$

Possibly taking a subsequence, we may consider (17) as a limit. For any $\varepsilon > 0$, we define

$$\Omega_\varepsilon := \{x \in \Omega : |\nabla \mathbf{h}| + |\nabla w| + |\nabla^2 w| < 1/\varepsilon, \ J > \varepsilon\}.$$

Clearly, $|\Omega \setminus \Omega_\varepsilon| \to 0$ as $\varepsilon \to 0$. By virtue of (H1), it follows that

$$\begin{aligned}
I[h_k, w_k] &\geq \int_{\Omega_\varepsilon} \Gamma(\nabla^2 w, \nabla \mathbf{h} + \mathbf{e}_3 \otimes \nabla w, J) dx \\
&+ \int_{\Omega_\varepsilon} \Gamma_K(\nabla^2 w, \nabla \mathbf{h} + \mathbf{e}_3 \otimes \nabla w, J) \cdot (\nabla^2 w_k - \nabla^2 w) dx \\
&+ \int_{\Omega_\varepsilon} \Gamma_F(\nabla^2 w, \nabla \mathbf{h} + \mathbf{e}_3 \otimes \nabla w, J) \cdot (\nabla \mathbf{h}_k - \nabla \mathbf{h}) dx \\
&+ \int_{\Omega_\varepsilon} \Gamma_F(\nabla^2 w, \nabla \mathbf{h} + \mathbf{e}_3 \otimes \nabla w, J) \cdot [\mathbf{e}_3 \otimes (\nabla w_k - \nabla w)] dx \\
&+ \int_{\Omega_\varepsilon} \Gamma_J(\nabla^2 w, \nabla \mathbf{h} + \mathbf{e}_3 \otimes \nabla w, J) \cdot (J_k - J) dx,
\end{aligned} \qquad (18)$$

where $\Gamma_\rho$ denotes partial differentiation with respect to the arguments $\rho = \mathbf{K}, \mathbf{F}, J$. On taking the limit as $k \to \infty$ in (18), we find that all integrals on the right side of the inequality vanish except for the first, by virtue of weak convergence. Thus, there remains

$$\lim_{k \to \infty} I[w_k, h_k] \geq \int_\Omega \mathcal{X}_{\Omega_\varepsilon} \Gamma(\nabla^2 w, \nabla \mathbf{h} + \mathbf{e}_3 \otimes \nabla w, J) dx, \qquad (19)$$

where $\mathcal{X}_{[\cdot]}$ denotes the characteristic function. In view of (8) and (H1), we see that $\Gamma(\cdot)$ is non-negative. Consequently, (17) follows from (19) in the limit as $\varepsilon \to 0$ in (19) via the monotone convergence theorem (and by the weak continuity of $\varphi$).

To finish the proof, we combine the above results to conclude

$$E[\mathbf{h}_*, w_*] \leq \liminf_{k \to \infty} E[\mathbf{h}_k, w_k] = \inf_\mathcal{A} E[\mathbf{h}, w] \leq E[\mathbf{h}_o, w_o] < \infty,$$

with $(\mathbf{h}_*, w_*) \in \mathcal{A}$. $\square$

## 4. Concluding Remarks

As already discussed in Section 2, (7) is consistent with the small-bending approximation (6). Once (6) is adopted, we have no control over the outward unit normal field $\mathbf{n}$. Consequently, the exact expression (4) is not attainable, and (H3) represents a consistent approximation of (10).

If we drop (H3) altogether, then it is a straightforward exercise to take a rigorous first variation to obtain the weak form of the Euler-Lagrange equilibrium equations. The caveat is that such solutions could possibly exhibit local orientation reversal in that case. In the presence of (H3), we can take an out-of-plane variation without disturbing the possibility that $\det \nabla \mathbf{h}_* = 0$ (implying $W \sim +\infty$) on a set of measure zero in $\Omega$. To see this, we first define the set

$$\Omega_n = \{x \in \Omega : \sup[|W_K(\nabla^2 w_*, \nabla \mathbf{h}_* + \mathbf{e}_3 \otimes \nabla w_*)|$$
$$+ |W_F(\nabla^2 w_*, \nabla \mathbf{h}_* + \mathbf{e}_3 \otimes \nabla w_*)|] \leq n, \ \det \nabla \mathbf{h}_* \geq 1/n\},$$

for any $n \in \mathbb{N}$. Next, define $E_n[\mathbf{h}, w]$ as before in (11) but with the integrals now restricted to $\Omega_n$. Then well-known arguments, e.g., [4], enable a rigorous first-variation condition, viz., we can take the limit $\lim_{t \to 0}\{E_n[w_* + t\eta, \mathbf{h}_*] - E_n[w_*, \mathbf{h}_*]\}/t,$ leading to

$$\int_\Omega \mathcal{X}_{\Omega_n}\{W_K(\nabla^2 w_*, \nabla \mathbf{h}_* + \mathbf{e}_3 \otimes \nabla w_*) \cdot \nabla^2 \eta$$
$$+ ([W_F(\nabla^2 w_*, \nabla \mathbf{h}_* + \mathbf{e}_3 \otimes \nabla w_*)]^T \mathbf{e}_3) \cdot \nabla \eta\} dx = 0, \ \forall \eta \in W_\Gamma^{2,p}(\Omega),$$

where $n$ is arbitrarily large. This represents a weak form of the out-of-plane equilibrium equation.

The existence of an energy minimizer for any of the theories employed in [2], [7], [13] and [14] does not follow from our approach; in each case, the membrane energy does not satisfy hypothesis (H1). In addition, our existence theorem does not extend to a theory like that proposed in [13] where the exact product (5) (with $\mathbf{n} = [\mathbf{f}_{,1} \times \mathbf{f}_{,2}]/J$) is employed instead of (6). We mention that an existence theorem for a Cosserat shell theory featuring, inter alia, membrane properties like those considered in the present work but based on (4), is provided in [9]. Certain sub-determinants involving products of components of both the deformation gradient and the gradient of a director field (a special case of which is the unit-normal field) are allowed in the theory, based on polyconvexity assumptions. However, the relationship of the latter with the model in [13] is unclear.

## Acknowledgements

This work was supported in part by the National Science Foundation through grant DMS-2006586, which is gratefully acknowledged. I thank Gokul Nair for useful comments on an earlier version of the work.